\documentclass{article}

\usepackage[final]{Arxiv_v1}

\usepackage[utf8]{inputenc} 
\usepackage[T1]{fontenc}    
\usepackage{hyperref}       
\usepackage{url}            
\usepackage{booktabs}       
\usepackage{amsfonts}       
\usepackage{nicefrac}       
\usepackage{microtype}      

\usepackage[pdftex]{graphicx}
\usepackage{amssymb,amsfonts,amsmath}
\usepackage{color}

\definecolor{white}{rgb}{1,1,1}
\definecolor{blue}{rgb}{0.1,0.1,1}
\definecolor{cyan}{rgb}{0.55,0.6,0.8}  
\definecolor{red}{rgb}{0.8,0.2,0.2}  
\definecolor{green}{rgb}{0,0.8,0} 
\definecolor{purple}{rgb}{0.6,0,1}
\definecolor{orange}{rgb}{1,0.6,0.3}

 \newcommand{\blue}{\color{black}} \newcommand{\cyan}{\color{black}} \newcommand{\red}{\color{black}}  \newcommand{\green}{\color{black}}

\title{Generalized Phase Representation of \\ Integrate-and-Fire Models}

\author{
Dongsung~Huh
}

\begin{document}

\maketitle

\begin{abstract}

The quadratic integrate-and-fire (QIF)  model captures the normal form bifurcation dynamics of Type-I neurons found in cortex. 
Remarkably, 
this model is known to have a dual equivalent representation in terms of phase, called the $\theta$-model, 
which has advantages for numerical simulations and analysis over the QIF model.
%
Here, I investigate the nature of the dual representation and derive the general phase model expression for all integrate-and-fire models. 
Moreover, I show the condition for which the phase models 
exhibit discontinuous onset firing rate, 
the hallmark 
of Type-II spiking neurons.

\end{abstract}

\section{Introduction}

Integrate-and-fire models 
provide the simplest dynamic descriptions of biological neurons 
that produce brief pulse signals to communicate,  
called spikes.  
They {\green  abstract out}  
the complex biological mechanism of spike generation 
with a simple threshold-reset process
\cite{lapicque_recherches_1907,izhikevich2007dynamical}, 
which allows 
efficient analysis and simulation of 
spiking  neural network models.

Generally, an integrate-and-fire model is described by 
a continuous subthreshold dynamics and a discontinuous threshold-reset process:
\begin{align}
	\label{eq:dyn_IF}
	\dot{x}  & =  f(x) +  I       & ( \text{for } x< x_+) \\
	x  & \to x_-                     & ( \text{for }  x = x_+) \nonumber 
\end{align}
where $x$ represents the  neuron's internal state, 
$I(t)$ is the net-input  ({\it i.e.} including rheobase)
and 
$x_+, ~ x_-$ are the threshold and the reset, respectively. 
For example,  
$f(x) =  |x| $,   $x_+ = 0 $,  $x_- =-1$
describes the leaky integrate-and-fire model
\cite{stein1965theoretical}.


%
An interesting example is  the quadratic integrate-and-fire (QIF) model, 
defined by $f(x)=x^2$ 
and the threshold and reset at infinity: $x_\pm = \pm \infty$
\cite{latham2000intrinsic}.
Remarkably, this model 
is known to have 
a dual representation in terms of phase, 
called the {\it $\theta$-model}
\cite{ermentrout1996type,ermentrout2008ermentrout}: 
\begin{align*}
	\dot{\theta} & =  (1-I) \sin(\theta)^2 +I  & ( \text{for } \theta< \pi/2) \\
	\theta  & \to -\pi/2                     & ( \text{for }  \theta = \pi/2)
\end{align*}
%
which transforms to QIF through the mapping $x = \tan(\theta) $.  
This phase representation 
maps the infinite real line to a unit circle $S^1$ 
{\blue through {\cyan one-point} compactification \cite{alexandroff1924metrisation,wiki:Compactification}
(see section 4.2),} 
providing 
the canonical model for Type-I cortical neurons that exhibit SNIC bifurcations (saddle-node bifurcation on an invariant circle)  \cite{ermentrout2008ermentrout}. 
%
Note that the phase representation 
removes the cumbersome infinities involved with 
the QIF model,
as well as effectively removing the {\green discontinuous} threshold-reset process,
due to the the periodicity of the model: $\sin(\theta)^2 = \sin(\theta + \pi)^2$.
%
Thus, the resulting phase dynamics is smooth and bounded, 
which is advantageous for numerical simulations and analysis of spiking neural network dynamics
\cite{monteforte2010dynamical,lajoie2013chaos,lajoie2014structured,novikov2016robustness},
as well as optimization \cite{huh2017gradient,mckennoch2009spike}.

However, the mathematical implications and the generalized form of the dual representation  
has not been investigated for general integrate-and-fire models.
In this manuscript, I investigate the nature of the dual representation 
and derive the following generalized phase dynamics  for all integrate-and-fire models:
\begin{align}
	\label{eq:dyn_IF_dual}
	\dot{y}  & =  (1-I) \, g(y) + I       & ( \text{for } y< y_+) \\
	y  & \to y_-                     & ( \text{for }  y = y_+) \nonumber 
\end{align}
which transforms to eq~\eqref{eq:dyn_IF} through a corresponding mapping function:  
$$x = h(y).$$ 

\section{ Derivation} 

{\red Here,} I introduce the idea of representative phase and dual ODEs 
to derive the dual form eq~\eqref{eq:dyn_IF_dual} of the integrate-and-fire dynamics.

\subsection{Representative phase}

The nonlinear function $f(x)$ of the integrate-and-fire models eq~\eqref{eq:dyn_IF}
is assumed to have a well-defined minimum value.
Without loss of generality,  a normalized form can be imposed 
to set  the minimum of $f(x)$ at the origin, such that 
$\min f(x) = f(0)=0$, 
%
with the threshold   $x_+ \geq 0 $, and the reset $x_- < 0$.
%
Also, $f(x)$ is assumed to diverge to infinity in the limit $x \to \pm\infty$, but not for finite $x$:  $f(x)<\infty, ~ \forall |x|<\infty$.
%
%
Note that  a positive  input $I>0$  drives the 
dynamics  to  fire spikes repeatedly (oscillatory regime), 
whereas a negative input $I<0$  ceases such on-going activities  (excitable regime). 

{\red Now,} consider 
the integrate-and-fire dynamics eq~\eqref{eq:dyn_IF} 
under a constant unit input $I=1$
as a representative case,
and define $y$ as the phase of this representative dynamics: 
{\it i.e.} 
\begin{align}
	\label{eq:dyn_2}
	 \frac{dx}{dy}  & =  \tilde{f}(x) , \\
	\text{where} ~~~~~~~~~ \tilde{f}(x) & \equiv f(x) +1   . ~~~~~~~~~~~~~~ \nonumber
\end{align}
Note that 
$y$ is defined as a {\it time-like} variable, and  $\tilde{f}(x) \geq 1$ for $\forall x$.

By definition, the phase dynamics under is  $\dot{y}=1$ for the constant  unit input $I=1$. 
Generalizing the phase dynamics to arbitrary time-varying input $I(t)$
 requires finding the dual ODE of eq~\eqref{eq:dyn_2}.

\subsection{Dual ODEs}

Dual ODEs are defined here as 
a pair of ordinary differential equations (ODEs) 
%
%
\begin{align*}
	 \frac{dx}{dy}  & =  \tilde{f}(x)     \\ 
%
	 \frac{dy}{dx}  & =  \tilde{g}(y)      
\end{align*}
%
that describe an identical relationship between $x$ and $y$. 
Such pairs can always be found for  monotonic relationships:
{\it i.e.} $\tilde{f}(x), \tilde{g}(y)  > 0$. 
In integral forms, the relationship is expressed as 
\begin{align}
	\label{eq:forward_map_dyn}
	x  & =  \int_0^y \frac{du}{\tilde{g}(u)}  \equiv h(y)  \\
	\label{eq:inverse_map_dyn}
	y &  =  \int_0^x \frac{du}{\tilde{f}(u)}   \equiv h^{-1}(x)  .
\end{align}
That is, the integral of $1/\tilde{g}$ is the inverse function of the integral of $1/\tilde{f}$.
Note that the condition 
$h(0)=0$
is imposed at the origin without loss of generality.

\subsection{Legendre transform}

The problem of finding the dual of an ODE is closely related to Legendre transform.
Legendre transform of a convex function $H(x)$ is 
\cite{zia2009making}
\begin{align}
	\label{eq:Legendre}
	H^*(y)  = \min_x \,   y x   - H(x) , 
\end{align}
where $y$ is  the dual  variable of $x$. 
%
The mapping between $x$ and $y$ is described by the first order derivatives
%
\begin{align*}
	y & = H'(x) 	~  \equiv h^{-1}(x) \\ 
	x & = {H^*}'(y)     \equiv  h(y) , 
\end{align*}
which  
can be identified with eq~(\ref{eq:forward_map_dyn},\ref{eq:inverse_map_dyn}) to yield   
\begin{align*}
	H''(x) & = {1}/{\tilde{f}(x)}       \\
	{H^*}''(y) & = {1}/{\tilde{g}(y)} .
\end{align*}
This result links Legendre transform 
to the dual ODE problem. 
Moreover, these second order derivatives are well known to be reciprocal to each other:
\begin{align}
	\label{eq:reciprocity_fg} 
	H''(x)    {H^*}''(y) 
	=  \frac{1}{\tilde{f}(x) \, \tilde{g}(y)}
	= 1.
\end{align}
%

\subsection {Generalized phase dynamics}

%
For arbitrary time-varying input $I(t)$, the phase dynamics generalizes to the following form
\begin{align}
	\label{eq:phase_dyn_temp}
	& & & & & &  & &  & & \dot{y} = \frac{dx}{dt} \frac{dy}{dx} & = ( f(x) +  I ) \tilde{g}( y) ,  ~~~   \\ 
	\text{where} ~~~~~~~~~~~~~~~~~~~~~~  & &  & &   & & & &  & & 
	  \frac{dy}{dx}  &  =  \tilde{g}(y)  & & & &  & & & & & &  ~~~~~~~~~~~~~~~~~~~~~~~~~~~~~~~~~~  \nonumber   
\end{align}
is the dual ODE of eq~\eqref{eq:dyn_2}, as determined by  Legendre transform. 
%
%

According to the reciprocity condition eq~\eqref{eq:reciprocity_fg}
and $f(x) = \tilde{f}(x) -1$,
equation~\eqref{eq:phase_dyn_temp} 
simplifies to 
\begin{align*}
	\dot{y} 	& = 1 +   (  I - 1 )\, \tilde{g}( y)  \\[2pt] 
			& =  (1 - I) \, g(y) + I,  
\end{align*}
which derives the generalized phase dynamics  eq~\eqref{eq:dyn_IF_dual}.
Here, $g(y) \equiv 1 - \tilde{g}(y)$ 
is defined to describe the intrinsic phase dynamics under the null input, $I=0$. 


\section{Examples}

\begin{figure} 
  \centering
   \includegraphics[width=0.99\textwidth]{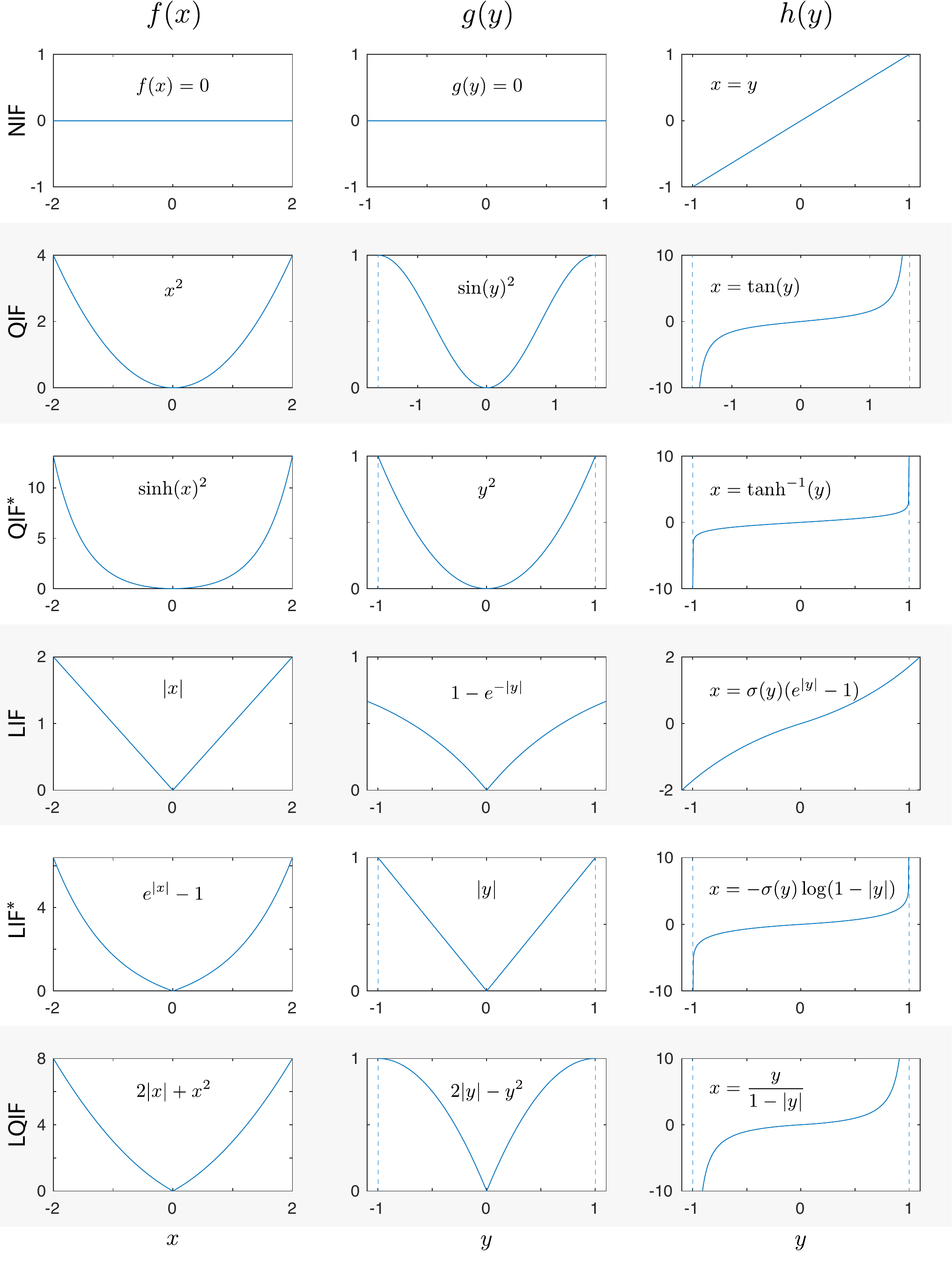}
	  \caption{Dynamic functions $f(x)$ of various integrate-and-fire models, the corresponding dual functions $g(y)$ of the phase dynamics, and the mapping functions $x=h(y)$ are shown.
	Dashed lines indicate where $h(y)$ diverge:  $ \lim_{x \to \pm \infty} h^{-1}(x)$.
	$\sigma(\cdot)$ is the sign function. 
	}
	\label{fig1}
\end{figure}

The above result generalizes 
the dual relationship 
between the QIF model and the $\theta$-model 
to all  integrate-and-fire models {\cyan of the form} 
eq~\eqref{eq:dyn_IF}. 
Some examples are shown here.  
See Figure~\ref{fig1} and Table~\ref{table1}.

\begin{itemize}

\item {\bf NIF }  Non-leaky integrate-and-fire  model:  Self-dual.
%
\begin{align*}
	~
	\dot{x} &  =  I , ~~~~ &
	\dot{y} & = I   ~~~~~~~~~~~~~~~~~~~~~~~~~~~~
\end{align*}

\item {\bf QIF }  Quadratic integrate-and-fire  model. 
%
\begin{align*}
	 ~~~~~
	\dot{x} &  =  x^2 + I , &
	\dot{y} & =  (1-I) \sin(y)^2 +I 
	~~
\end{align*}

\item {\bf QIF$^*$ }  Quadratic phase dynamics model. 
%
%
\begin{align*}
	 ~~~~~~~~~~~~~~
	\dot{x} &  =  \sinh(x)^2 + I, & 
	\dot{y} & =   (1-I) y^2 + I
	 ~~~~~~~~~~~~~~~~~~~~
\end{align*}

\item {\bf LIF } Linear (or Leaky) integrate-and-fire  model (symmetrized).
%
\begin{align*}
	 ~~~~~~~~~
	\dot{x} &  =  |x| + I ,  ~~~~ &
	\dot{ y} &  = (1-I) (1-e^{ - |y|}) + I 
\end{align*}

\item {\bf LIF$^*$ }  Linear phase dynamics model (symmetrized). 
%
%
\begin{align*}
	 ~~~~~~~~~~
	\dot{x} &  =  e^{|x|}-1 + I , &
	\dot{ y} &  = (1-I)  |y| + I 
	 ~~~~~~~~~~~~~~~
\end{align*}

\item {\bf LQIF } Linear-Quadratic integrate-and-fire  model  (symmetrized). 
%
\begin{align*}
	 ~~~~~~~~~~~~
	\dot{x} &  =    2|x| + x^2  + I ,  &
	\dot{ y} & =  (1-I) (2|y| - y^2) + I 
	~~~~
\end{align*}
%

%

\end{itemize}

\begin{table}
  \begin{tabular}{ | c || c | c | c | c | c | c |}
	\hline 
			    &  NIF & QIF & QIF$^*$  & LIF & LIF$^*$  & LQIF 
	\\ \hline \hline 
	$\tilde{f}(x)$  & $1$ & $1+x^2$    & $\cosh(x)^2$ 
			         & $1+|x|$ & $e^{|x|}$  & $(1+|x|)^2$       
	\\ [6pt]
	$\tilde{g}(y)$ & $1$ & $\cos(y)^2$ & $1-y^2$  
	                        & $\displaystyle e^{ - |y|}$ & $1-|y|$   & $(1-|y|)^2 $ 
	\\ [6pt]
	$h(y)$ & $y$ & $\tan(y)$      & $\tanh^{-1}(y) $      
	                        & $ \sigma(y)( e^{|y|} - 1)$  & $ - \sigma(y) \log (1 -  |y|)$        & $\displaystyle\frac{y}{1 -  |y| } $
	\\ [8pt]
	\hline 
  \end{tabular}
  \\ [4pt]
  \caption{Representative dynamic functions $\tilde{f}(x)$ of various integrate-and-fire models under unit input $I=0$, the corresponding dual functions $\tilde{g}(y)$ of the phase dynamics, and the mapping  functions $x=h(y)$ are shown. $\sigma(\cdot)$ is the sign function. }
	\label{table1}
\end{table}

\pagebreak

\section{Properties of the phase model}

\subsection{Phase dynamics function $g(y)$}

The reciprocity condition eq~\eqref{eq:reciprocity_fg} yields 
the following relationships between the dynamics functions of the dual representations
\begin{align}
	g(y) = \frac{f(x)}{1+f(x)}  
\end{align}
%
and between their derivatives
\begin{align}
	\frac{d g}{dy}  = \frac{df/dx}{1 + f}  =  \frac{d}{d x} \log (1 + f),
\end{align}
where $x=h(y)$ is implied.


Since ${f}(x) \geq 0 $, 
$g(y)$  is bounded by $0 \leq {g}(y) \leq 1$, 
%
%
%
with the minimum  at the origin, $f(0)=g(0)=0$.
In the limit $x\to \pm\infty$,  $g(y)$ approaches 1 as $f(x)$ diverges to infinity.

The derivatives  ${dg}/{dy} \approx {d f}/{d x} $  are similar near the origin, 
since $\log (1 + f(x)) \approx f(x)$ for $ |x| \ll 1$.
%
In the limit $x\to \pm\infty$,   $dg/dy$ 
{\red remains} finite only if $f(x)$ is an exponentially growing function  ({\it e.g.} QIF$^*$ and LIF$^*$).  
Otherwise, the derivative of $ \log (1 + f(x))$ approaches zero in the limit 
({\it e.g.}  polynomial $f(x)$: QIF, LIF, and LQIF).

\subsection{Compatification} 

The phase $y = h^{-1}(x)$ describes the time it takes for the representative dynamics $ \dot{x}=\tilde{f}(x)$  to reach $x$ from the origin.
%
Therefore, 
for 
model dynamics that diverges to infinity within finite time,
the phase $y_{\pm}$ that corresponds to the threshold/reset $x_\pm$ at infinity 
is finite
({\it e.g.} QIF, QIF$^*$,  LIF$^*$, LQIF).
Thus, for these fast diverging models,  the phase representation {\green shrinks} the infinite real line to an open interval $(y_-,y_+)$. 
%
The open interval can then be topologically {\green morphed} 
into a  circle $S^1$ by bringing the ends of the interval together and adding a single "point at infinity", a process called one-point compactification \cite{alexandroff1924metrisation,wiki:Compactification}.

The compactification of phase representations 
resolves the 
inconvenience of 
simulating diverging dynamics over the infinite domain. 
Moreover, 
it removes the discontinuity of the threshold-reset process
by equating the threshold phase with the reset phase, 
at which the dynamics is identical: 
$g(y_-) = g(y_+)=1$.  

\begin{figure} 
  \centering
   \includegraphics[width=0.69\textwidth]{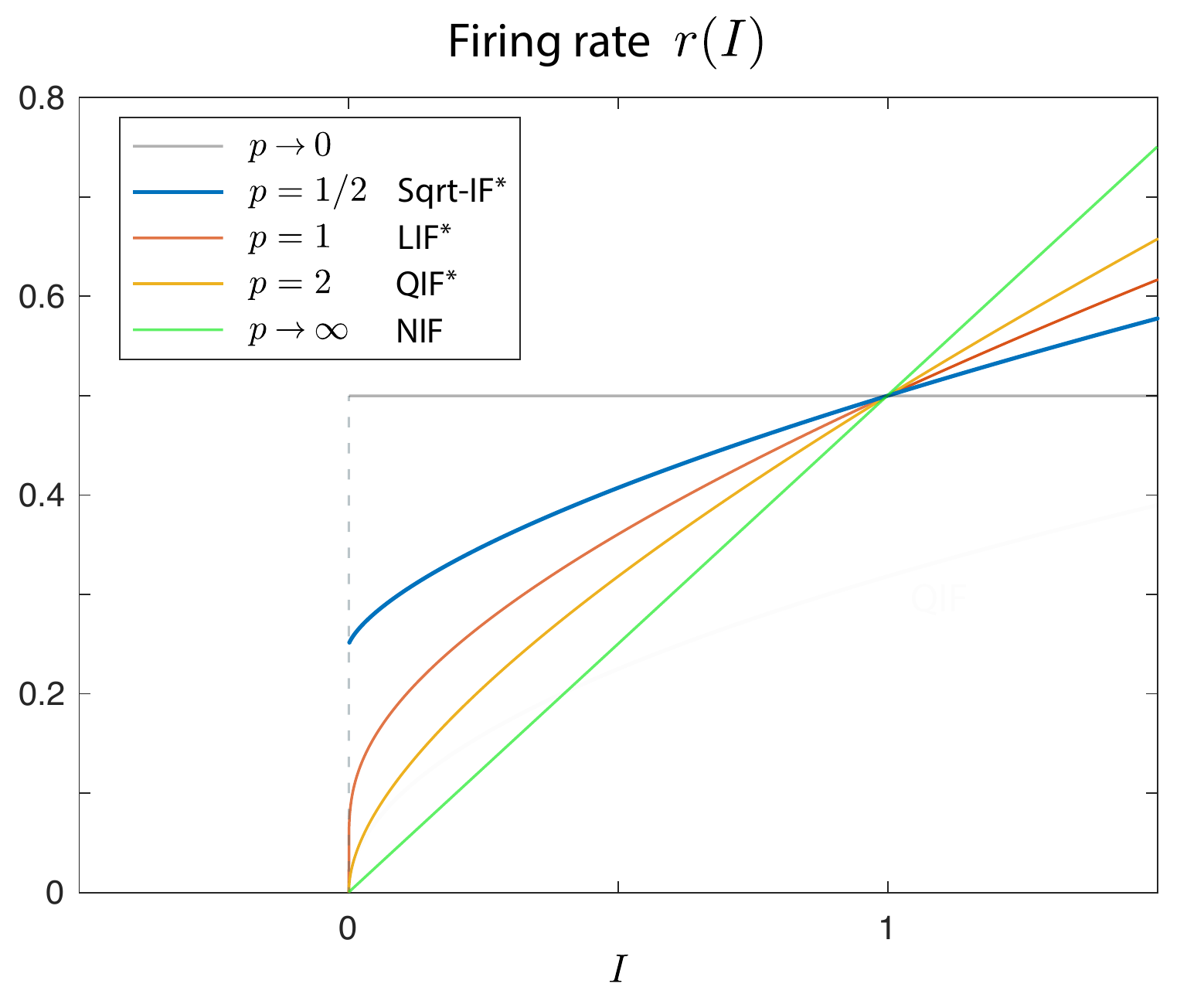} 
	  \caption{
			Firing rate of various phase dynamics models of the form $g(y) = |y|^p$.
			The limit $p\to\infty$ corresponds to the NIF model. 
			For $p\geq 1$, the onset firing rate is $r(0)=0$.
			For $p<1$, the models exhibit non-zero onset firing rate, 
			with discontiuous F-I cuves at $I=0$. 
			}
	\label{fig:FI_curve}
\end{figure} 

\begin{table}[b]
  \begin{tabular}{ | c || c | c | c | c | c | }
	\hline 
			    &  NIF & QIF$^*$  &LIF$^*$ & Sqrt-IF$^*$  &  LQIF  
	\\ \hline \hline 
	${g}(y)$ &  $0$  &  $y^2$ &  $|y|$    &  $\sqrt{|y|}$   & $2|y| - y^2 $ 
	\\ [6pt]
	$r(I)$  	&  $\displaystyle \frac I 2   $
			&  $\displaystyle \frac 1 2 \frac{\sqrt{ (I - 1) I } }{\tanh^{-1}(\sqrt{ (I - 1) / I})} $ 
			&  $\displaystyle \frac 1 2   \frac{I -1}{ \log I}$
			&  $\displaystyle \frac 1 4  \frac {(I - 1)^2} {1+I (\log I - 1) }  $
			&  $\displaystyle \frac 1 2 \frac{\sqrt {  I - 1}}{ \tan^{-1}(\sqrt{ I - 1}  ) } $ 
	\\ [12pt]
	\hline 
  \end{tabular}
  \\ [4pt]
  \caption{The firing rate  $r(I)$ 
of various phase dynamics models. $I\geq 0$. $y_\pm= \pm 1$.
%
}
  \label{table:FI}
\end{table}

\subsection{Firing rate} 

For positive constant input $I > 0$, integrate-and-fire models exhibit periodic spiking activities.
The spiking period can be calculated as 
$$ T (I) 
= \int_{x_-}^{x_+} \frac {d x}{   f(x) + I } 
= \int_{y_-}^{y_+} \frac {d y}{   (1 - I)  g(y) + I } .
$$
%
Note that the spiking period for unit input is always $T(1) = y_+ - y_-$  for all models. 

The inverse of period is called the firing rate: $r(I) \equiv 1/T(I)$. 
For negative constant input $I<0$, the firing rate is $r(I) =0$. 
The firing rate functions are {\green shown} in Table~\ref{table:FI} for various integrate-and-fire models. 

Figure~\ref{fig:FI_curve} shows 
the firing rate plotted against the input, called the F-I curves.
In many models, the firing rate continuously decreases to zero with the input, 
the defining characteristics of Type-I cortical neurons. 
%
Remarkably, however, 
the {\it square-root} phase dynamics model (Sqrt-IF$^*$), defined by $g(y)=\sqrt{|y|}$, exhibits non-zero firing rate at onset  $I=0$, 
%
the hallmark of Type-II neurons  \cite{skinner2013moving}.

More generally, the onset firing rate of monomial phase dynamics models 
%
	$ g(y)  = |y|^p $  
%
is
%
\begin{align*}
	 & & \lim_{I \to 0^+} r(I)   
		& =  \frac{1-p}{2}  & (\text{for } p<1) \\
	& & & = 0  & (\text{for } p\geq 1) 
\end{align*}
(with $y_{\pm}=\pm 1$), 
where the limit is taken from the positive side.
Thus, the F-I curves of monomial phase models of $p<1$  exhibit 
{Type-II-\it like}  discontinuities  at onset $I=0$. 


The onset characteristics 
 generalizes to other models 
whose phase dynamics function approximates 
$g(y) \approx |y|^p$ near the origin. 
%
For example,  
the integrate-and-fire model with  $f( x)   =  e^{\sqrt{2 |x|}} - 1$  exhibits non-zero onset firing rate%
\footnote{Its firing rate is $r(I)  =   (1 - I)/2\,\text{Li}_2(1 - I) $ with $r(0) =   {3}/{\pi^2}$. $\text{Li}_n(\cdot)$ is the poly-logarithmic function.},
since its phase dynamics is $g(y) \approx \sqrt{2|y|}$ for small $|y| \ll 1$.
{\red Also, 
QIF$^*$'s firing rate is approximated by QIF's rate $r(I)=\sqrt{I}/\pi$  near the onset,
and
LIF$^*$'s firing rate  is approximated by LIF's rate 
$r(I)=[2\log (1+1/I)]^{-1}$ near the onset  ($x_\pm=\pm1$).
}

\section{Summary}

I introduced the formal definition for the phase representation of integrate-and-fire models
and derived the generalized expression for the phase dynamics models. 
The compactified phase representation for integrate-and-fire models would facilitate efficient  simulations and analysis of large scale spiking neural network dynamics.  

{\red
Analysis of the phase representation revealed a new class of integrate-and-fire models that blur the distinction between Type-I and Type-II neurons \cite{skinner2013moving}. 
Despite the discontinuous onset of F-I curves, 
however,
these models 
do not truly qualify as Type-II,
since their  phase response curves are strictly positive (not shown here), 
and they exhibit SNIC bifurcations rather than Hopf or saddle-node-off-limit-cycle bifurcations.}
{\blue Remarkably, however, these models can also exhibit  
supercritical Hopf bifurcation with sustained subthreshold oscillations when paired with an adaptation process (not shown here), which require further investigation.}

\subsubsection*{Acknowledgments}
The author would like to thank Bard Ermentrout for insightful discussion.

\bibliographystyle{unsrt}

\end{document}